\documentclass[11pt]{article}
\usepackage{amsmath,amsthm,amssymb}
\usepackage{graphicx}
\usepackage{booktabs}
\usepackage{enumerate}
\usepackage{color}
\usepackage{titlesec}
\usepackage{wrapfig}
\titleformat{\section}
    [block]
    {\bfseries\large}
    {\thesection}
    {0.5em}
    {}
\titlespacing{\section}
    {0pt}
    {18pt}
    {6pt}

\newtheorem{Th}{Theorem}[section] 
\newtheorem{Lem}[Th]{Lemma} 
\newtheorem{Prop}[Th]{Proposition} 
 
\newtheorem{Def}[Th]{Definition} 
 
\newtheorem{Rem}[Th]{Remark} 
 
\renewcommand{\thesection}{\arabic{section}}

\newcommand{\argmin}{\mathop{\rm argmin}\limits}

\def\R{{\mathbb R}}

\def\cC{{\cal C}} 
 
\def\cE{{\cal E}}
\def\cF{{\cal F}}
 
\def\cH{{\cal H}}

\def\DS{\displaystyle} 
\def\1/2{\frac{1}{2}}

\def\DIV{\mbox{\rm div}}

\def\ov{\overline} 
 
\def\vep{\varepsilon}

\def\GammaD{{\Gamma_{\rm\scriptsize D}}}
\def\GammaN{{\Gamma_{\rm\scriptsize N}}}
\def\Gammaz{{\Gamma_{\rm\scriptsize N}^0}}
\def\Gammao{{\Gamma_{\rm\scriptsize N}^1}}

\makeatletter 
\long\def\@makefntext#1{\parindent 1em\noindent 
\@hangfrom{\hbox to 1.8em{\hss$^{\@thefnmark}$}}#1}
\makeatother

\begin{document}
\begin{center}
  {\large\bf What is the physical origin of the gradient flow structure\\
    of variational fracture models?
  }
\end{center}

\begin{center}
  Masato Kimura$^*$, Takeshi Takaishi$^{**}$, Yoshimi Tanaka$^{***}$\\~\\
\begin{tabular}{l}
*) Kanazawa University, mkimura@se.kanazawa-u.ac.jp\\
**) Musashino University, taketaka@musashino-u.ac.jp\\
***) Kanazawa Gakuin University, yoshimi-t@kanazawa-gu.ac.jp
\end{tabular}
\end{center}

\vspace{1cm}
\begin{abstract}
  We investigate a physical characterization of the gradient flow
  structure of variational fracture models for brittle materials: a Griffith-type fracture model
  and an irreversible fracture phase field model.
  We derive the Griffith-type
  fracture model by assuming that the fracture energy in Griffith's theory
  is an increasing
  function of the crack tip velocity. Such a velocity dependence
  of the fracture energy is typically observed in polymers.
  We also prove an energy dissipation identity of the Griffith-type fracture model,
  in other words, its gradient flow structure. On the other hand,
  the irreversible fracture phase field model is derived
  as a unidirectional gradient flow of a regularized total energy
  with a small time relaxation parameter based on the variational
  fracture theory by Francfort and Marigo (1998)
  and a mathematical space regularization proposed by Ambrosio and Tortorelli (1992).
  We have considered the time relaxation parameter a mathematical approximation parameter,
  which we should choose as small as possible. In this research, however,
  we reveal the physical origin of the gradient flow structure
  of the fracture phase field model and show that the small time relaxation
  parameter is characterized as the rate of velocity dependence of the fracture energy.
  It is verified by comparing the energy dissipation properties
  of those two models and by analyzing a traveling wave solution
  of the irreversible fracture phase field model.
\end{abstract}


\section{Introduction}\label{Sec:intro} 
\setcounter{equation}{0}
This paper considers variational fracture models for quasi-static crack propagation
in a brittle material, especially a variant of the Griffith-type fracture model
and an irreversible fracture phase field model.
We also discuss their energy dissipation
properties and the physical characterization of a small
time relaxation parameter in the variational fracture model,
that is, $\alpha >0$ in \eqref{fpfmeq0} below.

Bourdin et al. \cite{B-F-M00} and Karma et al. \cite{K-L-K01}
initiated the phase field approach to model fracture phenomena.
Then, it is widely used for numerical studies of
the dynamics of fracture under complex geometries and conditions in 2D or 3D,
such as fractures in thermoelasticity \cite{BKM2010, BMMS2014, KTANT2021, Alfat22},
or viscoelasticity \cite{KTANT2021,TT2020},
crack nucleation \cite{TLBMM2018,KBFL20},
and cracking phenomena with other physical and chemical effects
\cite{KTANT2021,MBGL13,CBY19}.

\begin{wrapfigure}{r}{0.50\textwidth}
    \centering
\includegraphics[width=6.5cm]{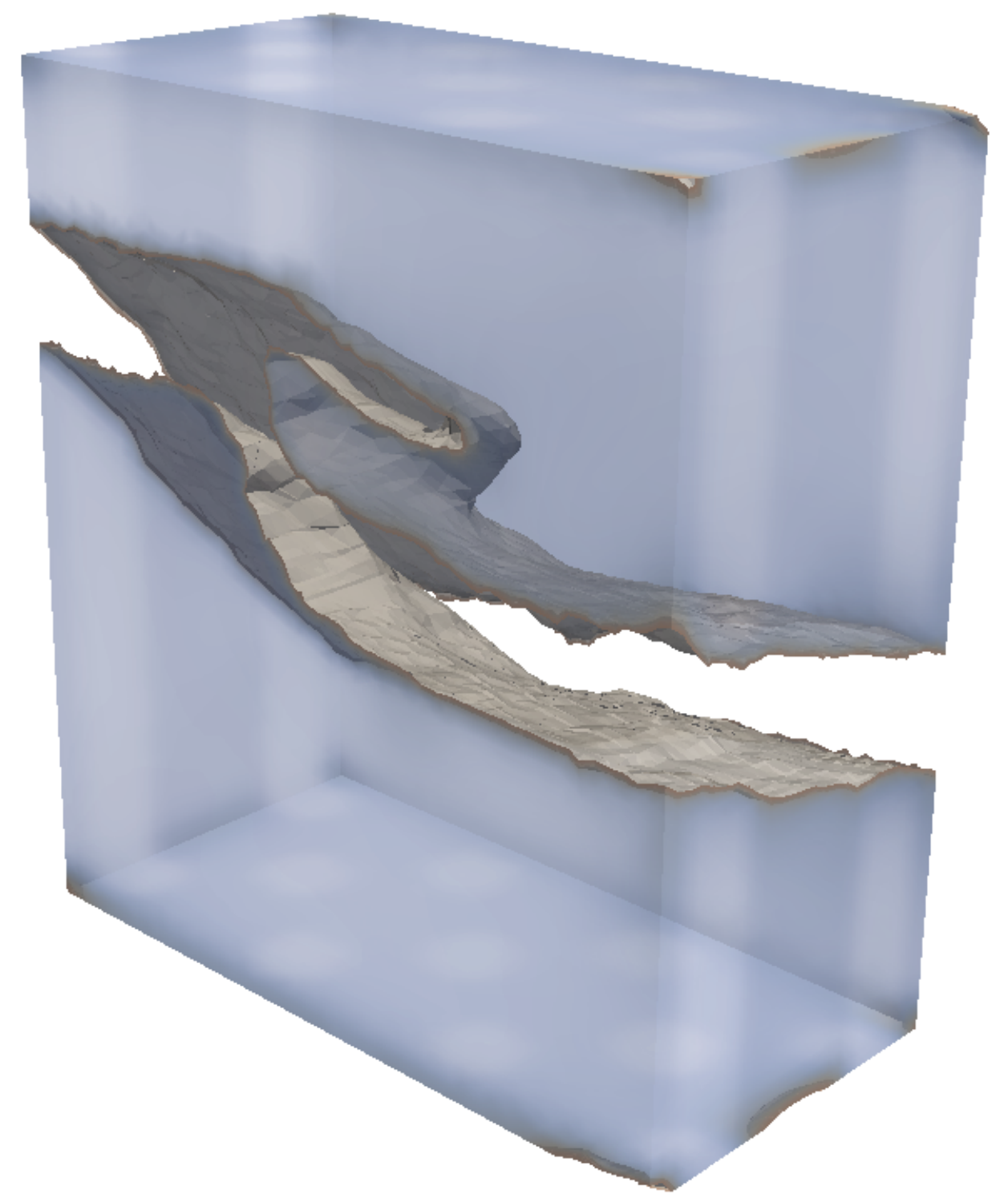}
\caption{An example of fracture simulation by F-PFM in 3D.}
    \label{fig:fpfm}
\end{wrapfigure}

The phase field model is a diffused interface approach to the crack problem,
i.e., instead of describing the crack as a sharp boundary,
a smooth phase field variable (damage field variable) is
introduced over the material region.
The following phase field model for fracture phenomena
(which is denoted by F-PFM in this paper)
was proposed in \cite{KTANT2021, T-K09}:
\begin{align}\label{fpfmeq0}
\begin{cases}
-\DIV \left((1-z)^2 \sigma [u] \right) = f(t),\\
\DS{\alpha \frac{\partial z}{\partial t}=
\left( \vep \; \DIV (G_c \nabla z ) 
       - \frac{G_c}{\vep} z + \sigma[u]:e[u] (1 - z) 
       \right)_+},
\end{cases}
\end{align}
where $u(x,t)\in\R^d$ ($d=2,3$) denotes a displacement
and $z(x,t)\in [0,1]$ denotes a phase field variable for
the crack position as $z\approx 1$ for the cracked region
and $z\approx 0$ for the undamaged region. The phase field $z$ is
often called a damage variable.
The parameters $\alpha$ and $\vep$
are small positive real numbers related to regularizations in time and space,
respectively.
As a crack can not be healed itself,
we take the positive part $(~)_+$ of the right-hand
side of the second equation,
where $(a)_+ = \max (a,0)$.
The use of the  positive part function guarantees the irreversibility of the crack propagation:
$\frac{\partial z}{\partial t}\ge 0$.
Figure~\ref{fig:fpfm} shows an example of a finite element simulation of a
complex fracture geometry by the F-PFM in 3D.
See more detail in Section~\ref{Subsec:fpfm-edi} and also \cite{KTANT2021}.

As shown in \cite{KTANT2021, Alfat22, TT2020, T-K09}, the F-PFM successfully modeled
various fracture phenomena with energy consistency.
The F-PFM includes two artificial small
positive parameters $\vep$ and $\alpha$, which relate to 
the space regularization and the time relaxation, respectively.
Roughly speaking, the crack tip singularity of the stress field
is regularized by $\vep$, and the ``sudden jump'' singularity (see Figure~\ref{fig:Gclen})
of the crack propagation is regularized by $\alpha$.
These regularizations enable us to get a stable numerical
crack propagation.
However, the physical characterization of these small parameters
have yet to be well studied.

This paper aims to clarify the physical characterization of
the parameter $\alpha>0$ in the F-PFM. As shown in \cite{KTANT2021},
the F-PFM satisfies an energy dissipation identity \eqref{dissipation},
and $\alpha$ becomes a coefficient of the dissipation term.

On the other hand,
forming the process zone near the crack tip/edge
causes such energy dissipation,
and it is experimentally observed as a velocity dependence of the
fracture energy (see Section~\ref{Subsec:Vdependent} for details).
To clarify the connection between the energy dissipation and the velocity dependence of the fracture energy, we consider a Griffith-type fracture model.  Then, we reveal that the velocity-dependent fracture energy causes energy dissipation
The obtained energy dissipation identity represents its gradient flow structure,
which resembles the one of the F-PFM.
Through such mathematical evidence, we systematically explain
that the parameter $\alpha$ in the F-PFM has a clear physical meaning
as the rate of velocity dependence of the fracture energy.

The outline of this paper is as follows.
Section~\ref{Sec:VFT} briefly reviews
the energy dissipation identities
in the classical Griffith theory and
the variational fracture theory when the crack path is prescribed.
Then, in Section~\ref{Sec:cpm-vdfg}, we consider Griffith's crack propagation
model with the velocity-dependent fracture energy and prove that
it can be described as a well-posed initial value problem of an ODE, and
satisfies a natural energy dissipation identity.
In Section~\ref{Subsec:fpfm-edi}, we will see that the gradient flow structure
of the F-PFM implies an energy dissipation identity that resembles
one of the ODE models in Section~\ref{Sec:cpm-vdfg}.
Furthermore, we investigate the regularized fracture energy
of the F-PFM by considering a traveling wave solution in Section~\ref{Subsec:TS}.
Section~\ref{Subsec:trqs} 4.3 also mentions a physical interpretation
of another time relaxation parameter initially introduced by
\cite{T-K09}. Finally, we will give concluding remarks and open questions in the last section.

\section{Quasi-static variational fracture theory}\label{Sec:VFT} 
\setcounter{equation}{0}

\subsection{Crack problem in linear elasticity}\label{Subsec:crack-prob}
We first consider a crack problem in static linear elasticity.
We omit details of notation and mathematical assumptions here
and refer to Section~2 of \cite{KTANT2021} for more precise
definitions and mathematical settings. In this paper, for simplicity, we often
abbreviate the space variable $x$, e.g., $u(t)$ means $u(t)=u(x,t)$ or $u(t)=u(\cdot,t)$.

\begin{wrapfigure}{r}{0.50\textwidth}
  \centering
  \includegraphics[width=0.49\textwidth]{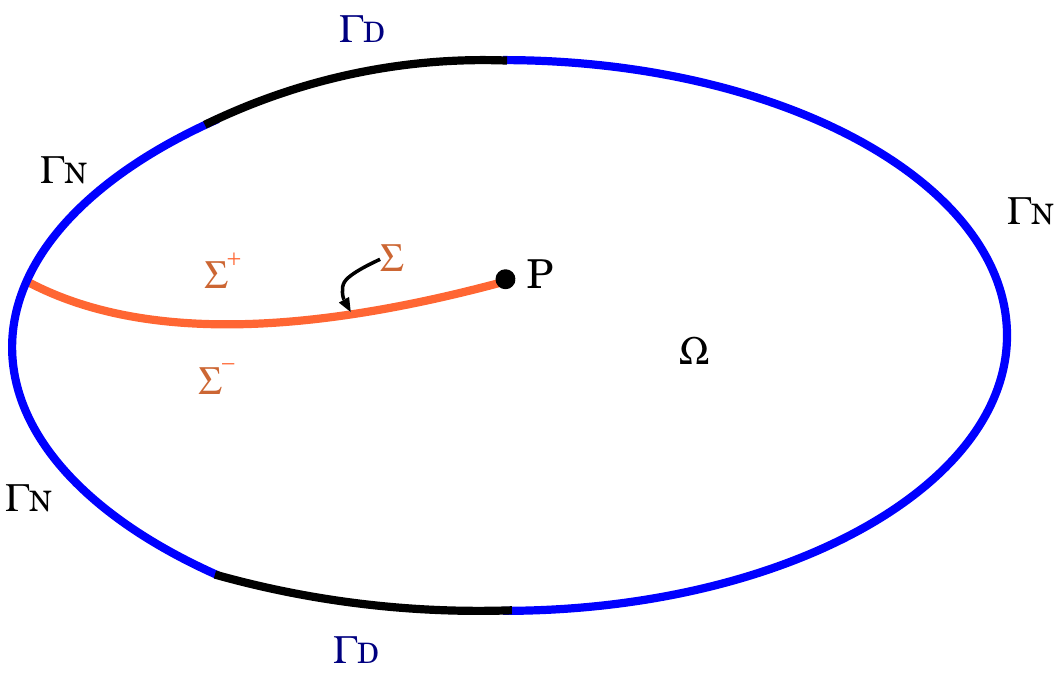} 
 \caption{A cracked domain $\Omega\setminus\Sigma$ with a Dirichlet boundary $\GammaD$
   and a Neumann boundary $\GammaN$.}\label{domain}
\end{wrapfigure}

Let $\Omega$ be a bounded domain in $\R^d$ $(d=2,3)$,
which represents an elastic body.
We suppose a crack  $\Sigma\in \cC_0$ in $\Omega$,
where $\cC_0$ is an admissible set of cracks (see Appendix~\ref{Sec:cecp}, and also
Section~2 of \cite{AKA}).
We denote the length/area (for $d=2$/$d=3$) of the crack $\Sigma$ by $|\Sigma |$.
We consider the following boundary value problem of linear elasticity in
the cracked domain $\Omega\setminus \Sigma$:

\begin{equation}\label{elas}
\left\{
\begin{array}{rcll}
- \mbox{div } \sigma [u] & = & f(t)\quad & \mbox{in }\Omega\setminus\Sigma,\\[5pt]
\sigma [u] \nu & = & q(t) & \mbox{on } \GammaN ,\\[5pt]
\sigma [u] \nu & = & 0 & \mbox{on } \Sigma^\pm ,\\[5pt]
u & = & g(t) & \mbox{on } \GammaD ,
\end{array}
\right .
\end{equation}
Under suitable assumptions, there exists a weak solution to \eqref{elas},
and we denote it by $u(t;\Sigma):\Omega\setminus\Sigma \to \R^d$.
It is known that the following
variational principle gives $u(t;\Sigma)$:
\begin{align}\label{ut}
u(t;\Sigma)=\argmin_{v\in V(g(t);\Sigma)} E_{el}(t,v;\Sigma),
\end{align}
where $V(g;\Sigma):=\{ v\in H^1(\Omega\setminus \Sigma;\R^d);(v-g)|_{\GammaD}=0\}$
for $g\in H^1(\Omega;\R^d)$, and
\begin{align}\label{Eel}
  E_{el}(t,v;\Sigma):=
\frac{1}{2} \int_{\Omega\setminus\Sigma} \sigma[v]:e[v]\,dx
-\int_\Omega f(t)\cdot v\,dx-\int_{\GammaN} q(t)\cdot v\,ds
\end{align}
represents the elastic energy of a displacement $v\in H^1(\Omega\setminus\Sigma;\R^d)$
including the body and surface forces at time $t$.
Then, the elastic energy in the cracked body $\Omega\setminus\Sigma$
at time $t$ is given as
\begin{align}\label{min}
E_{el}^*(t;\Sigma):=
\min_{v\in V(g(t);\Sigma)} E_{el}(t,v;\Sigma)
=E_{el}(t,u(t;\Sigma);\Sigma),
\end{align}
and it is known (e.g.\cite{F-M98}) that
\begin{align}\label{mono}
  E_{el}^*(t;\Sigma)\ge E_{el}^*(t;\tilde{\Sigma})\quad\mbox{holds, if}\quad
  \Sigma\subset \tilde{\Sigma}\in \cC_0.
\end{align}

When we fix the crack $\Sigma$ and the given loads $(g(t),f(t),q(t))$ change in time smoothly,
the following energy conservation property holds:
\begin{align}\label{dtE}
  \frac{d}{dt}E_{el}^*(t;\Sigma)=\dot{F}(t,u(t;\Sigma);\Sigma),
\end{align}
where
\begin{align}\label{dotF}
\dot{F}(t,v;\Sigma):=\int_{\GammaD}\frac{\partial g}{\partial t}(t)\cdot (\sigma[v]\nu) \,ds
-\int_{\Omega\setminus\Sigma}\frac{\partial f}{\partial t}(t)\cdot v\,dx
-\int_{\GammaN}\frac{\partial q}{\partial t}(t)\cdot v\,ds.
\end{align}
The three terms on the right-hand side of \eqref{dotF}
represent the rates of energy injection for a displacement $v$
through the boundary displacement $g(t)$, the body force $f(t)$,
and the surface traction $q(t)$, respectively.
Using the integration
by parts formula under suitable regularity assumptions,
we can derive 
the energy identity \eqref{dtE}.

\subsection{Energy profile and energy release rate along a given crack path}
\label{Subsec:ep}
In the pioneering work by A. A. Griffith \cite{Grif20},
he constructed an energetic fracture theory under
the assumption that a crack path is given and that crack evolution is continuous in time.
Please refer to Appendix~\ref{Sec:cecp} for the precise definitions of a crack path
$\{\Sigma_p(\ell)\}_{\ell_0 \leq \ell \leq \ell_1}$
and a crack evolution $\{\Sigma(t)\}_{t_0 \leq t \leq t_1}$.

For a given crack path $\{\Sigma_p(\ell)\}_{\ell_0\leq \ell \leq \ell_1}\subset \cC_0$,
which is parametrized by $\ell = |\Sigma_p(\ell)|$,
we define
$E(\ell,t) := E_{el}^*(t;\Sigma_p(\ell))$ for
$\ell \in [\ell_0, \ell_1 ]$,
and refer to the function $\ell\mapsto E(\ell, t)$ 
as an energy profile.
If the energy profile $E(\ell,t)$ is of $C^1$-class in $\ell$,
then 
$G(\ell, t):= -\frac{\partial E}{\partial \ell}(\ell,t)$
is called an energy release rate per unit length/area of the crack evolution.
From \eqref{mono}, it follows that $E(\ell,t)$ is nonincreasing in $\ell$ and 
$G(\ell, t)\geq 0$ holds.
From \eqref{dtE}, we also have
\begin{align}\label{dtEF}
  \frac{\partial E}{\partial t}(\ell,t)=\dot{F}(t,u(t;\Sigma_p(\ell));\Sigma_p(\ell)).
\end{align}

\subsection{Griffith theory}\label{Subsec:G-model}
According to \cite{B-F12, B-F-M08, C-F-M09},
the classical  Griffith theory is summarized as follows. We suppose that $\{\Sigma(t)\}_{t_0 \leq t \leq t_1}$ is a smooth, continuous crack evolution in $\Omega$ under a given boundary condition $g(t)$, and $\{\Sigma_p(\ell)\}_{\ell_0 \leq \ell \leq \ell_1}$ is the corresponding
crack path (Proposition~\ref{prop:crackPathProp}).
Then, there exists $G_c > 0$ such that
$L(t):=|\Sigma(t)|$ satisfies the following conditions:
\begin{equation}\label{Griffith}
  \begin{cases}
    L'(t) \geq 0\quad \text{(Irreversibility)} \\
    G(L(t), t) \leq G_c\quad \text{(Griffith's Criterion)}\\
    L'(t)\left(G_c -G(L(t), t) \right)  = 0\quad \text{(Energy Conservation)} 
  \end{cases}
\end{equation}
for $t \in [t_0, t_1]$,
where $G_c$ is a material property called a fracture energy (or a critical energy release rate).
The third condition of \eqref{Griffith} represents the conservation of a
total energy:
\begin{align}\label{Etot}
  E_{tot}^*(t;\Sigma):=E_{el}^*(t;\Sigma)+G_c|\Sigma|.
\end{align}
From \eqref{dtEF} and $G(\ell, t):= -\frac{\partial E}{\partial \ell}(\ell,t)$, we have
\begin{align*}
\frac{d}{dt}E_{tot}^*(t;\Sigma(t))
=\frac{d}{dt}(E(L(t),t)+G_cL(t))
=L'(t)(G_c-G(L(t),t))+\dot{F}(t,u(t;\Sigma(t));\Sigma(t)).
\end{align*}
This implies the following energy conservation law:
\begin{align*}
\frac{d}{dt}E_{tot}^*(t;\Sigma(t))=\dot{F}(t,u(t;\Sigma(t));\Sigma(t)),
\end{align*}
provided the third condition of \eqref{Griffith} holds.

\section{Crack propagation model with velocity-dependent fracture energy}
\label{Sec:cpm-vdfg}
\setcounter{equation}{0}

\subsection{Velocity-dependent fracture energy}\label{Subsec:Vdependent}
Many experiments \cite{And05,Evans1974,DFI1985,Owen1998} on metals,
ceramics, and polymers have revealed that the measured fracture energy
(or, equivalently, critical stress intensity factor) depends on crack 
velocity.
We denote the crack velocity-dependent fracture energy by $G^*_c(V)$,
where $V\ge 0$ is the crack tip velocity in 2D and
the normal component (i.e., normal to the crack edge)
of the crack edge velocity in 3D.

The physical origin of the $V$-dependence is the formation
of the so-called process zone around the crack tip \cite{And05}. 
A process zone has an intermediate spatial scale
(far larger than the atomic scale and far smaller than the specimen size),
and some dissipative processes occur there.
The size of the process zone and the intensity of the energy dissipation change with $V$,
and we can macroscopically measure those dependencies on $V$
as a $V$-dependence of the fracture energy.

\begin{figure}[htb]
\centering
\includegraphics[width=1.0\textwidth]{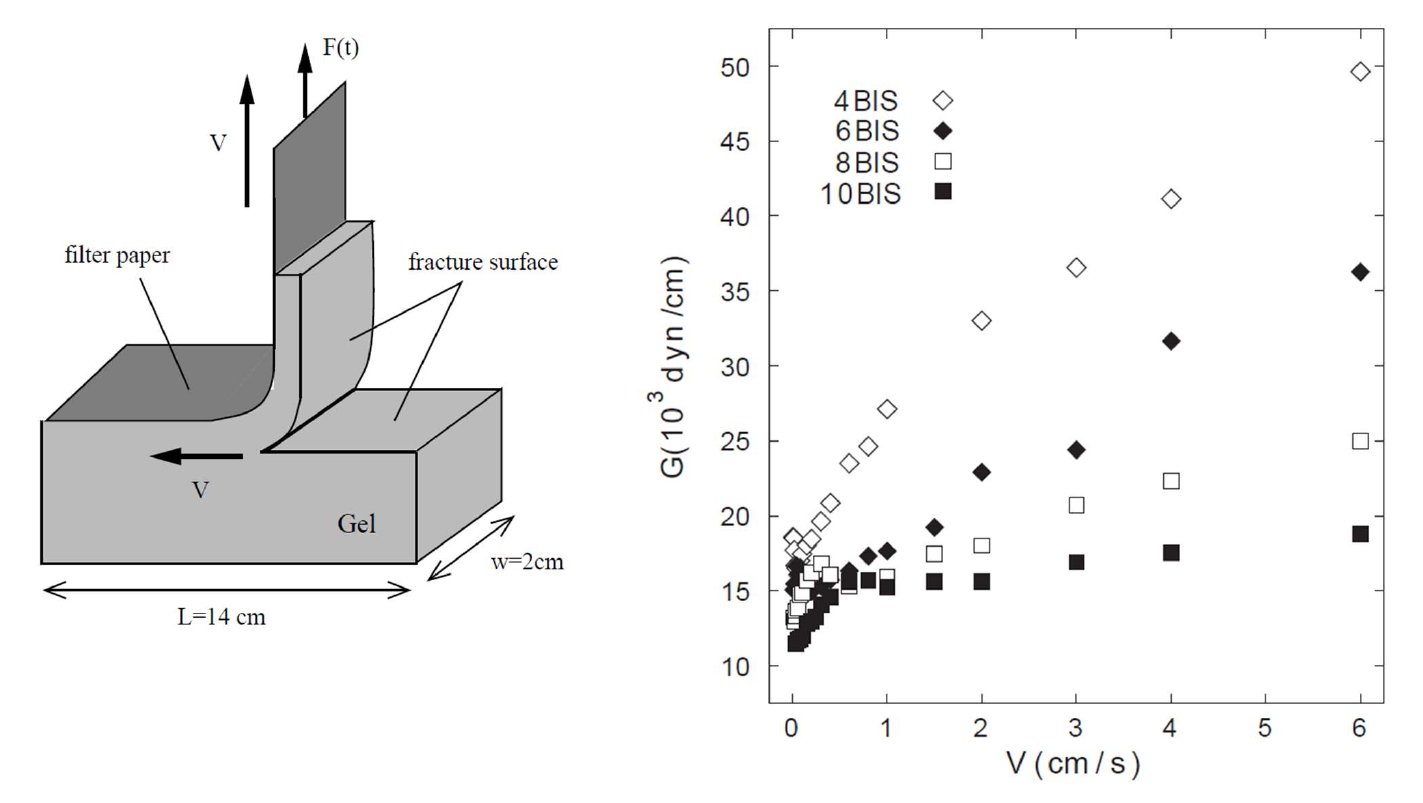}
\caption{Crack velocity ($V$) dependence of fracture energy ($G$) of chemically-crosslinked acrylamide hydrogels measured by a sort of tearing test (the left illustration), taken from \cite{TFM2000} with kind permission of The European Physical Journal: The difference in plot symbols represents the difference in crosslink density. As the crosslink density increases, $G(V)$ gets lower. The data of $G(V)$ show slightly upper convex behavior but almost linear for larger $V$.}
    \label{fig:gelG(v)}
\end{figure}

Usually, $G^*_c(V)$ increases with $V$
(the faster deformations cause the larger dissipations).
Especially, gel materials, crosslinked polymer networks swollen
with solvent, tend to show a simple, almost linearly increase behavior,
as seen in Fig. \ref{fig:gelG(v)} \cite{Baum2006,TFM2000}.
In the following lines, we assume $G^*_c(V)$ is a strictly increasing function of $V$.
However, it is experimentally possible
that $G^*_c(V)$ shows a negative slope or a drastic drop in a particular
$V$ region if the fracture mechanism qualitatively changes
in the $V$ region (e.g., brittle-ductile transition \cite{And05}).

\subsection{An ODE model and energy dissipation}\label{Subsec:ODEmodel} 
In this section, we set $d=2$ and suppose the crack has a single tip $\mbox{P}$,
as in Figure~\ref{domain}.
Then, $V=L'(t)$ denotes the crack propagation velocity.
We assume the following condition on $G_c^*(V)$:
\begin{align}\label{HGcV}
  \begin{cases}
  G_c^*(V)=G_c+\alpha^* (V)\quad (V\in [0,\infty)),\quad G_c>0,\\
    \mbox{$\alpha^*$ is a strictly increasing continuous function
      on $[0,\infty)$ with $\alpha^*(0)=0$.}
  \end{cases}
\end{align}

When the fracture energy depends on the crack tip velocity as
$G_c^*(V)$,
the Griffith model \eqref{Griffith} becomes
\begin{equation}\label{G2}
  \begin{cases}
    V\geq 0,\\
    G\leq G_c^*(V),\\
    V\left(G_c^*(V) -G\right)  = 0,
  \end{cases}
\end{equation}
where $G=G(L(t),t)$.
We call \eqref{G2} a Griffith-type fracture model with velocity-dependent fracture energy.
We have the following theorem.
\begin{Th}\label{P1}
Under the condition \eqref{HGcV}, 
the velocity-dependent fracture energy model \eqref{G2}
is equivalent to
\begin{align}\label{aV}
  \alpha^* (V)=(G-G_c)_+.
\end{align}
It is also equivalent to 
\begin{align}\label{V=beta}
  V=\beta^* (G-G_c),
\end{align}
where
\begin{align*}
  \beta^* (s):=
  \begin{cases}
    0&(s< 0)\\
    (\alpha^*)^{-1}(s)&(s\ge 0)
  \end{cases}
\end{align*}
\end{Th}
\proof
Under the condition \eqref{HGcV},
$\alpha^*(V)\ge 0$ holds if and only if $V\ge 0$ holds,
and $\alpha^*(V)= 0$ holds if and only if $V= 0$ holds.
So, \eqref{G2} is equivalent to
\begin{equation}\label{G2a}
  \begin{cases}
    \alpha^*(V)\ge 0,\\
    G_c^*(V)-G\ge 0,\\
    \alpha^*(V)\left(G_c^*(V)-G\right)  = 0.
  \end{cases}
\end{equation}
Then, applying \eqref{ab+} of Lemma~\ref{lemabc}, we obtain that
\eqref{G2a} is equivalent to
\begin{align*}
  \alpha^*(V)=(\alpha^*(V)-(G_c^*(V)-G))_+=(G-G_c)_+.
\end{align*}
The equivalency to the alternative form \eqref{V=beta} is quickly confirmed.
\qed

\begin{Rem}\label{wp}
{\rm
From Theorem~\ref{P1}, the Griffith-type model \eqref{G2} with initial crack length $\ell_0$
is equivalent to the following initial value problem of an ODE:
\begin{equation}\label{G2I2}
  \begin{cases}
    L'(t)=\beta^* (G(L(t),t)-G_c)\quad (t\ge t_0),\\
    L(t_0)=\ell_0.
  \end{cases}
\end{equation}
The function $G(\ell,t)$ is assumed to be continuous in $(\ell,t)$
and locally Lipschitz in $\ell$. If $\beta^*$ is also locally
  Lipschitz (for example, this is true if
  $\alpha^*\in C^1([0,\infty))$ and $(\alpha^*)'(V)>0$ for $V\ge 0$), then
    from the Cauchy-Lipschitz theorem,
    it follows that there exists a unique solution to \eqref{G2I2}
    locally in time.
}
\end{Rem}

\begin{figure}[htb]
    \centering
\includegraphics[width=10cm]{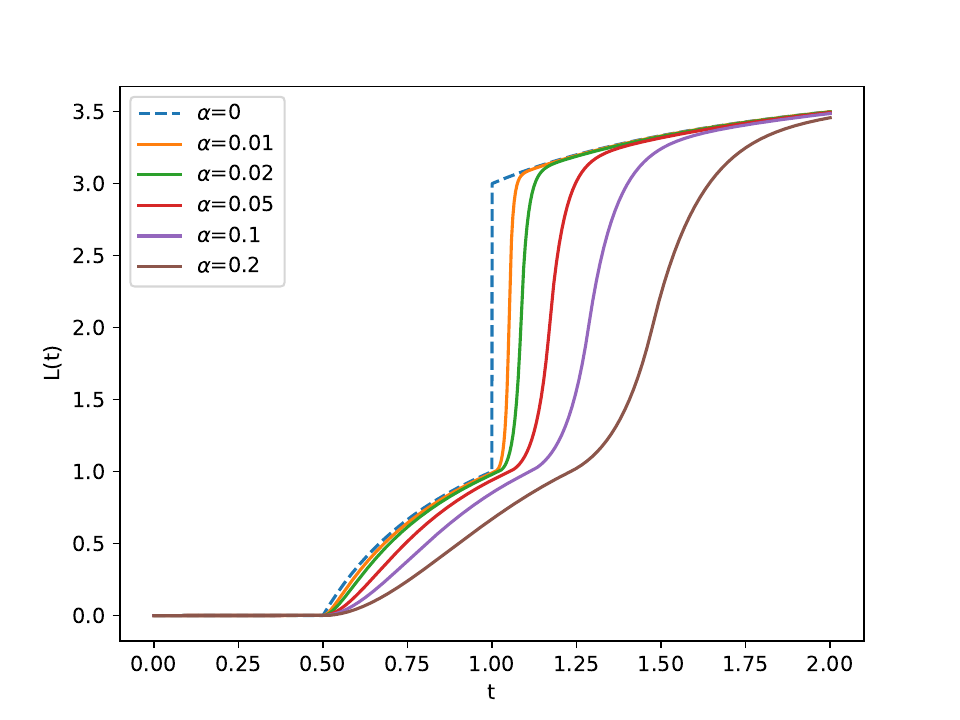}
\caption{Numerical solutions of \eqref{aODE} for $\alpha = 0.01\sim 0.2$
  with $G(\ell,t)=t(2-||l-1|-1|)$ and $G_c=1$. The broken line shows the limit
profile of $L(t)$ as $\alpha\to 0$, which has a sudden jump at $t=1$.}
    \label{fig:Gclen}
\end{figure}

\begin{Rem}\label{alphaV}
  {\rm
When $\alpha^*(V)$ is a linear function as $\alpha^*(V)=\alpha V$ with $\alpha >0$,
then $\beta^*(s)=\frac{1}{\alpha}(s)_+$ holds. In this case, \eqref{G2I2} becomes
\begin{equation}\label{aODE}
  \begin{cases}
    \alpha L'(t)=(G(L(t),t)-G_c)_+\quad (t\ge t_0),\\
    L(t_0)=\ell_0.
  \end{cases}
\end{equation}
In Figure~\ref{fig:Gclen}, we draw numerical solutions
of \eqref{aODE} for different $\alpha\in [0.01,0.2]$
with $G_c=1$ and an artificially given energy release rate function
$G(\ell,t):=t(2-||l-1|-1|)$. The broken line in the figure shows the limit
profile of $L(t)$ as $\alpha\to 0$, which has a sudden jump at $t=1$.
Such a sudden jump in the crack propagation is described in the framework of the
variational fracture theory by Francfort and Marigo \cite{F-M98,B-F-M08}.
However, the limit profile in the figure captures a slightly
different behavior from the original variational fracture theory \cite{F-M98}.
It corresponds to a localized Franctort-Marigo model \cite{AKTM15,AKA}.
We no longer discuss this issue in this paper, but
it will be discussed in our forthcoming paper intensively.
}
\end{Rem}

\begin{Th}[Energy dissipation identity]\label{th-edi}
Under the settings in Section~\ref{Subsec:G-model},
we suppose that the energy profile $E(\ell,t)$ satisfies
$E\in C^1([\ell_0,\ell_1]\times [t_0,t_1])$ and \eqref{dtEF},
and that $L\in C^1([t_0,t_1])$ be a solution of \eqref{G2I2} on $[t_0,t_1]$.
We define $\Sigma(t):=\Sigma_p(L(t))$ and $V(t):=L'(t)$.
Then, it satisfies the following energy dissipation identity:
\begin{align}\label{edi}
  \frac{d}{dt} E_{tot}^*(t;\Sigma(t))
  =-\alpha^*(V(t))V(t)+\dot{F}(t,u(t;\Sigma(t));\Sigma(t)).
\end{align}
In particular, when $\alpha^*(V)=\alpha V$, 
\begin{align*}
  \frac{d}{dt} E_{tot}^*(t;\Sigma(t))
  =-\alpha |V(t)|^2+\dot{F}(t,u(t;\Sigma(t));\Sigma(t)).
\end{align*}
\end{Th}
\proof
Since $E_{tot}^*(t;\Sigma(t))=E(L(t),t)+G_cL(t)$, we have
\begin{align*}
  \frac{d}{dt} E_{tot}^*(t;\Sigma(t))
  &=\frac{d}{dt}\left(E(L(t),t)+G_cL(t)\right)\\
  &=\left(\frac{\partial E}{\partial \ell}(L(t),t)+G_c\right)L'(t)
  +\frac{\partial E}{\partial t}(L(t),t)\\
  &=-\left(G(L(t),t)-G_c\right)V(t)
  +\dot{F}(t,u(t;\Sigma(t));\Sigma(t)).
\end{align*}
Hence, \eqref{edi} follows from 
$(G(L(t),t)-G_c)V(t)=\alpha^*(V(t))V(t)$, which we derive from
\eqref{aV}.
\qed

\subsection{Discussion}\label{subsec:discussion}
In this section, we studied the Griffith-type model \eqref{G2} with
a velocity-dependent fracture energy $G_c^*(V)=G_c+\alpha^*(V)$,
typically observed in polymer materials.
Theorem~\ref{P1} proved that the Griffith model with the velocity-dependent
fracture energy \eqref{G2} is equivalent to the ODE model \eqref{G2I2}.
In particular, when the fracture energy linearly depends on the velocity
$G_c^*(V)=G_c+\alpha V$, then it is written in the form:
$\alpha L'(t)=(G(L(t),t)-G_c)_+$ and satisfies the
energy dissipation identity:
$  \frac{d}{dt} E_{tot}^*(t;\Sigma(t))=-\alpha |V(t)|^2+\dot{F}$.
As we will see in the next section, the above
energy dissipation structure closely resembles the
one of the fracture phase field model.

\section{Irreversible fracture phase field model (F-PFM)}\label{Sec:FPFM}
\setcounter{equation}{0}

\subsection{F-PFM and energy dissipation identity}\label{Subsec:fpfm-edi}
\begin{figure}[h]
    \centering
    \includegraphics[clip,width=5cm]{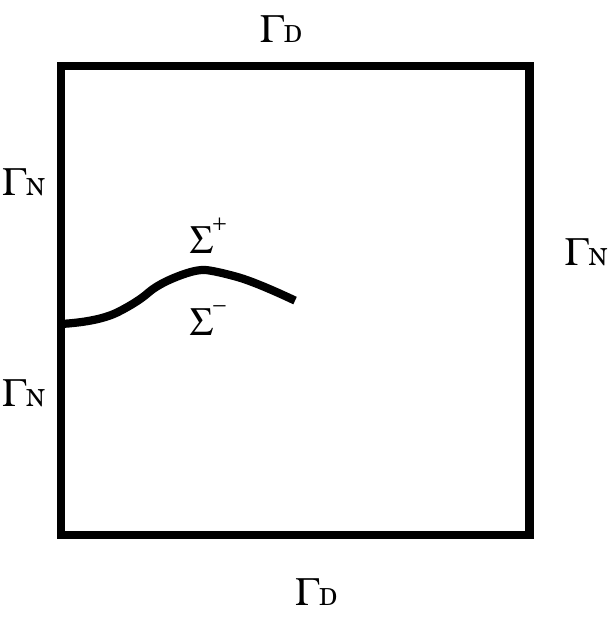}\hspace{10pt}
        \includegraphics[clip,width=6cm]{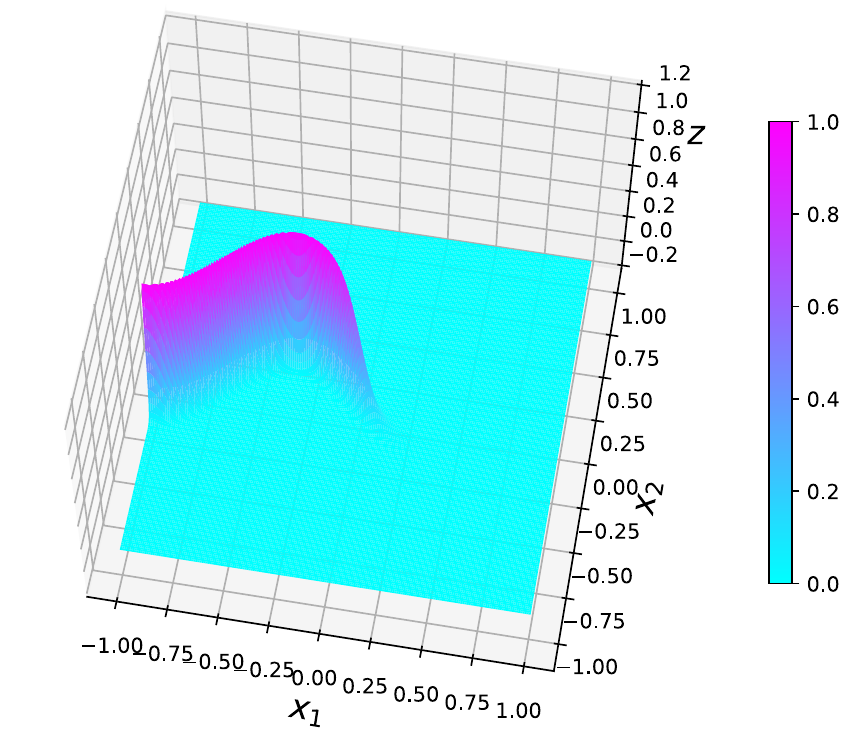}
        \caption{A crack $\Sigma$ in a two-dimensional
          rectangular domain $\Omega$ and a corresponding phase field variable $z(x)$
          are illustrated in the left and right figures, respectively.}
    \label{fig:crackz}
\end{figure}

This section briefly introduces an irreversible fracture phase field model (F-PFM)
based on \cite{KTANT2021, T-K09}.
We consider a smooth phase field function $z(x,t)$ 
to represent an approximate profile of the crack $\Sigma (t)$ (Fig.\ref{fig:crackz}).
We assume that $0\le z(x,t)\le 1$ and $z(x,t)\approx 1$
around crack $\Sigma(t)$, and that $z(x,t)\approx 0$ for the other region.
The phase field $z$ is also called a damage variable, representing a relative amount
of the accumulated damage in the elastic material. With the damage variable $z$,
$\tilde{C}:=(1-z)^2C$ gives the damaged elasticity tensor,
where $C$ denotes the original non-damaged elasticity tensor.

The F-PFM is described as the following initial and boundary value problem
of an elliptic-parabolic system of partial differential equations:
\begin{align}\label{fpfmeq}
\begin{cases}
-\DIV \left((1-z)^2 \sigma [u] \right) = f(t)\quad &\mbox{in}~\Omega\times [0,T],\\
\DS{\alpha \frac{\partial z}{\partial t}=
\left( \vep \; \DIV (G_c \nabla z ) 
       - \frac{G_c}{\vep} z + \sigma[u]:e[u] (1 - z) 
       \right)_+}\quad &\mbox{in}~\Omega\times (0,T],\\
u  =  g(t)\quad &\mbox{on}~\GammaD\times [0,T],\\
\sigma[u]\nu   =  q(t)\quad &\mbox{on}~\GammaN \times [0,T],\\
\DS{\frac{\partial z}{\partial \nu}  = 0}\quad &\mbox{on}~\Gamma\setminus\Gammao \times [0,T],\\
z  = 0\quad  &\mbox{on}~\Gammao\times [0,T],\\
z|_{t=0}  =  z^0\quad &\mbox{in}~\Omega .
\end{cases}
\end{align}
We suppose that $q(x,t)=0$ for $x\in \Gammaz\subset \Gamma_N$ and $t\in [0,T]$,
and set $\Gammao:=\GammaN\setminus \Gammaz$.
The second equation, a
nonlinear parabolic equation of $z$, describes the crack propagation.
The parameters $\alpha$ and $\vep$
are small positive real numbers related to regularizations in time and space,
respectively.
The positive part of the second equation's right-hand side
guarantees the crack propagation's irreversibility.

Instead of the elastic energy $E_{el}(t,v;\Sigma)$ of \eqref{Eel} and
the surface energy $G_c|\Sigma|$, we define the following
regularized elastic energy $\cE_{el}(t,v,z)$ and surface energy $\cE_s(z)$
applying the Ambrosio-Tortorelli approximation \cite{A-T92}:
\begin{align*}
&\cE_{el}(t,u,z):=\frac{1}{2} \int_\Omega (1-z)^2\sigma[u]:e[u]\,dx
-\int_\Omega f(t)\cdot u\,dx-\int_{\GammaN} q(t)\cdot u\,ds,\\
&\cE_s(z):=\frac{1}{2} \int_\Omega G_c
\left(
\vep |\nabla z|^2+
\frac{1}{\vep}z^2
\right)\,dx.
\end{align*}
We set $V(g):=\{ v\in H^1(\Omega;\R^d);(v-g)|_{\GammaD}=0\}$
for $g\in H^1(\Omega;\R^d)$.
Similarly to the case of $E_{el}$, we define
\begin{align*}
  &u(t,z):=\argmin_{v\in V(g(t))} \cE_{el}(t,v,z),\\
  &\cE_{el}^*(t,z):=\min_{v\in V(g(t))} \cE_{el}(t,v,z)
  = \cE_{el}(t,u(t,z), z),\\
  &\cE_{tot}^*(t,z):=\cE_{el}^*(t,z)+\cE_s(z).
\end{align*}
The F-PFM \eqref{fpfmeq} is derived
as a so-called irreversible gradient flow \cite{A-K19}
of $\cE_{tot}^*(t,z)$ with respect to $z$:
\begin{align*}
  \alpha \frac{\partial z}{\partial t}=\left(-\frac{\delta \cE_{tot}^*}{\delta z}\right)_+.
\end{align*}
In \cite{KTANT2021}, the following energy dissipation equality was shown.
If $(u(t),z(t))$ is a sufficiently smooth solution of \eqref{fpfmeq}, then it satisfies
\begin{align}\label{dissipation}
  \frac{d}{dt}\cE_{tot}^*(t,z(t))
  =-\alpha \int_\Omega \left|\frac{\partial z}{\partial t}\right|^2\,dx +\dot{\cF}(t,u(t),z(t)),
\end{align}
where
\begin{align*}
  \dot{\cF}(t,u,z):=
  \int_{\GammaD}\frac{\partial g}{\partial t}(t)\cdot ((1-z)^2\sigma[u]\nu) \,ds
-\int_{\Omega}\frac{\partial f}{\partial t}(t)\cdot u\,dx
-\int_{\GammaN}\frac{\partial q}{\partial t}(t)\cdot u\,ds.
\end{align*}

As shown above, the F-PFM is derived based on the
variational fracture theory
\cite{Grif20,F-M98},
the Ambrosio--Tortorelli regularization
\cite{A-T92}, and the unidirectional gradient flow
\cite{A-K19},
and it exhibits a natural energy dissipation property \eqref{dissipation} consequently.
In contrast to the other crack propagation models,
the F-PFM implicitly includes the crack path search and
enables us to treat complex crack patterns even in 3D (Figure~\ref{fig:fpfm}).

On the other hand, from the viewpoint of physics,
there are two open questions about the modeling of F-PFM.
One is the physical characterization of the damage variable $z$
and the spatial regularization parameter $\vep$. Here, $z$ and $\vep$ are
introduced in the mathematical regularization technique \cite{A-T92},
and their phyical substances have not been clarified yet.

The other open question is a physical characterization of the
time relaxation parameter $\alpha$, introduced in the
gradient flow. This paper aims to clarify the physical meaning
of the parameter $\alpha$ in the F-PFM. As we have discussed in Section~\ref{Sec:cpm-vdfg},
$\alpha$ in the ODE model \eqref{aODE} is characterized by the rate of
velocity dependence of the fracture energy $G_c^*(V)$: $\alpha =\frac{dG_c^*}{dV}$.
From the strong analogy between \eqref{aODE} and F-PFM, we expect
to characterize the parameter $\alpha$ in F-PFM similarly.
To strengthen this claim, we study the regularized fracture energy
of the F-PFM by considering a traveling wave solution in the following subsection.

\subsection{Traveling wave solution and velocity dependence of the fracture energy}\label{Subsec:TS} 
\begin{figure}[bht]
    \centering
\includegraphics[width=10cm]{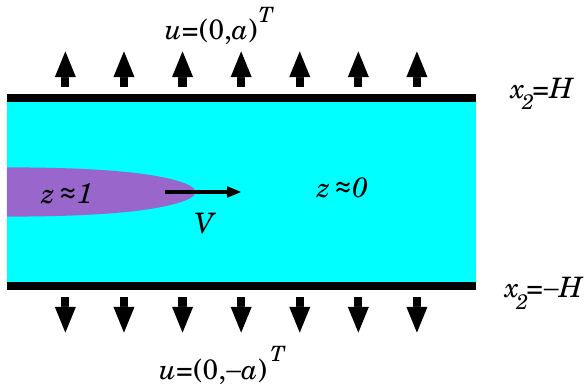}
\caption{A traveling wave solution of the F-PFM in 2D.}
    \label{fig:tws}
\end{figure}

In this section, we consider an infinite strip domain
$\Omega_H := \R\times (-H,H)\subset \R^2$ as shown in Figure~\ref{fig:tws},
and consider a traveling wave solution of the F-PFM in $\Omega_H$.
We set $\Gamma_H^\pm :=\{x=(x_1,x_2)^T\in\R^2;~x_2=\pm H\}$.
We consider the F-PFM in the strip domain $\Omega_H$ for $t\in\R$.
\begin{equation}\label{bvp}
\begin{cases}
-\DIV \left((1-z)^2 \sigma [u] \right) = 0\quad &\mbox{in}~\Omega_H\times \R,\\
\DS{\alpha \frac{\partial z}{\partial t}=
\left( G_c\left( \vep \Delta z - \frac{z}{\vep}\right)  + (1-z)\sigma[u]:e[u] \right)_+} 
     \quad &\mbox{in}~\Omega_H\times \R,\\
     u  = \left(\begin{array}{c}0\\\pm a\end{array}\right)
       \quad &\mbox{on}~\Gamma_H^\pm\times \R ,\\
\DS{\partial_2 z  = 0}\quad &\mbox{on}~\Gamma_H^\pm\times \R,\\
\end{cases}
\end{equation}
where $a>0$.

This geometry corresponds to the so-called ``pure sure geometry''
experimentally realized. Moreover, the fracture energy's velocity dependence
for rubbers and gels is measured in that experimental setting \cite{S-LP2023}.
The experimental system is a long rectangular plate $(-W,W)\times (-H,H)\times (0,b)$,
where $b>0$ is the thickness of this plate.
We realize a plane stress state in condition $b<<H<<W$.
An initial crack is made from the center of the left edge along the horizontal axis
($x_1$ axis).
We suppose that a pair of constant vertical boundary displacements of $\pm a$
at the upper and bottom edges of the system, i.e., $u_{2}(x_1, \pm H, t)=\pm a$, is applied.
The above boundary value problem \eqref{bvp} corresponds to the case $W=\infty$.

We suppose that there exists a traveling wave solution with a velocity $V>0$
in the direction of positive $x_1$, i.e.,
there exists $\bar{u}:\ov{\Omega_H}\to\R$ and $\bar{z}:\ov{\Omega_H} \to (0,1)$
and $V>0$ 
such that
\begin{align*}
u(x_1,x_2,t)=\bar{u}(x_1-Vt,x_2)=\bar{u}(\xi),\quad
z(x_1,x_2,t)=\bar{z}(x_1-Vt,x_2)=\bar{z}(\xi),
\end{align*}
where we set a moving coordinate $\xi:=(x_1-Vt,x_2)^T\in \ov{\Omega_H}$.
Then,
since $\partial_t z(x,t)=-V\partial_1\bar{z}(\xi)$, 
$(\bar{u},\bar{z}, V)$ is a solution of the following system:
\begin{equation}\label{tweq}
\begin{cases}
-\DIV \left((1-\bar{z})^2 \sigma [\bar{u}] \right) = 0\quad &\mbox{in}~\Omega_H,\\
\DS{-\alpha V\partial_1 \bar{z} =
G_c\left( \vep \Delta \bar{z} - \frac{\bar{z}}{\vep}\right)  + (1-\bar{z})w } 
     \quad &\mbox{in}~\Omega_H,\\
     \bar{u}  = \left(\begin{array}{c}0\\\pm a\end{array}\right)
       \quad &\mbox{on}~\Gamma_H^\pm,\\
\DS{\partial_2 \bar{z}  = 0}\quad &\mbox{on}~\Gamma_H^\pm,\\
\end{cases}
\end{equation}
where we have defined the density of elastic energy $w(\xi)$ by
$w(\xi):=\sigma [\bar{u}]:e[\bar{u}](\xi)$ for $\xi\in\ov{\Omega_H}$.
Additionally, we omitted the positive part of the second equation since
the expected profile of the traveling wave solution is $z_t >0$.

According to a number of our numerical experiments of the F-PFM,
we expect a traveling wave solution that corresponds to the constant-velocity
crack propagation, as shown in Figure~\ref{fig:tws}.
Since \eqref{tweq} is shift-invariant in the direction of $\xi_1$,
for a solution $(\bar{u}(\xi),\bar{z}(\xi), V)$,
$(\bar{u}(\xi_1-c,\xi_2),\bar{z}(\xi_1-c,\xi_2), V)$
is also a solution of \eqref{tweq} for any $c\in\R$.
We fix a solution $(\bar{u}(\xi),\bar{z}(\xi), V)$ in the following discussion.

For $t>0$ and a sufficiently small $\vep >0$,
the increment of the regularized crack length during
the time interval $(0,t)$ is given by
\begin{align*}
  L_\vep(t):=\1/2 \int_{\Omega_H}
  \left\{ \left(\vep |\nabla z(x,t)|^2+\frac{z(x,t)^2}{\vep}\right)
-\left(\vep |\nabla z(x,0)|^2+\frac{z(x,0)^2}{\vep}\right)\right\}
  \,dx .
\end{align*}
Since $V$ should coincide with $L_\vep'(t)$, it holds that
\begin{align}\label{LtV}
  L_\vep'(t)\approx V.
\end{align}
We also have
\begin{align*}
  L'_\vep(t)
  &= \int_{\Omega_H}  \left\{ \vep \nabla z(x,t)\cdot \nabla z_t(x,t)
  +\frac{1}{\vep}z(x,t)z_t(x,t)\right\}\,dx\\
  &= -\int_{\Omega_H}  \left( \vep \Delta z(x,t)-\frac{1}{\vep}z(x,t)\right)z_t(x,t)\,dx\\
  &= V\int_{\Omega_H}  \left( \vep \Delta \bar{z}-\frac{\bar{z}}{\vep}\right)
  \partial_1\bar{z}\,d\xi .  
\end{align*}
Similarly, 
the increment of the elastic energy during
the time interval $(0,t)$ is estimated by
\begin{align*}
  E_\vep(t):=\1/2 \int_{\Omega_H}
  \left\{ \left((1-z(x,t))^2\sigma[u(\cdot,t)]:e[u(\cdot,t)]\right)
-  \left((1-z(x,0))^2\sigma[u(\cdot,0)]:e[u(\cdot,0)]\right)
\right\}
\,dx ,
\end{align*}
and we have
\begin{align*}
  E'_\vep(t)
  &=
  -\int_{\Omega_H}
  ((1-z(x,t))\partial_tz(x,t)\sigma[u(\cdot,t)]:e[u(\cdot,t)]\,dx\\
&\quad + \int_{\Omega_H}(1-z(x,t))^2\sigma[u(\cdot,t)]:e[\partial_tu(\cdot,t)]
  \,dx\\
  &=
  V\int_{\Omega_H}
  ((1-\bar{z}(\xi))\partial_1\bar{z}(\xi)w(\xi)\,d\xi
  - \int_{\Omega_H}
  \DIV \left((1-z(x,t))^2\sigma[u(\cdot,t)]\right)\cdot \partial_tu(\cdot,t)
  \,dx\\
  &=
  V\int_{\Omega_H}
  ((1-\bar{z})w\partial_1\bar{z}\,d\xi\\
  &=
  -V\int_{\Omega_H}
  \left\{
  \alpha V\partial_1 \bar{z} 
  + G_c\left( \vep \Delta \bar{z} - \frac{\bar{z}}{\vep}\right)
  \right\}\partial_1\bar{z}(\xi)\,d\xi
\end{align*}

The increment of the total energy during
the time interval $(0,t)$ is given by $E_\vep(t)+G_cL_\vep(t)$ and
we have the energy dissipation identity:
\begin{align}\label{ediV}
  \frac{d}{dt}\left(E_\vep(t)+G_cL_\vep(t)\right)
  =-\alpha \beta V^2,
\end{align}
where
\begin{align*}
  \beta :=\int_{\Omega_H} |\partial_1 \bar{z}|^2\,d\xi >0.
\end{align*}

Since the fracture energy (the critical energy release rate) $G_c>0$ is defined by
the ratio of the released elastic energy per unit length
of the propagating crack, we consider an effective fracture energy $G_\vep$
for the traveling wave solution $\bar{z}(\xi)$ of the F-PFM:
\begin{align*}
  G_c^\vep :=-\frac{E'_\vep(t)}{V}.
\end{align*}
From the energy dissipation identity \eqref{ediV} and the approximation \eqref{LtV},
we obtain
\begin{align}\label{Gcep}
  G_c^\vep=G_c\frac{L_\vep '(t)}{V}+\alpha \beta V
  \approx G_c+\alpha \beta V.
\end{align}
The obtained formula \eqref{Gcep} suggests that the time relaxation parameter $\alpha$
in the F-PFM corresponds to the rate of velocity dependence of the regularized
fracture energy $\alpha\approx \frac{d G_c^\vep}{dV}$.

\subsection{Time relaxation for quasi-stationary elasticity}\label{Subsec:trqs}
Before concluding remarks, we briefly discuss the possible
modification of the F-PFM
on the quasi-stationary elasticity equation.
We consider a dynamic fracture model, replacing the first equation of \eqref{fpfmeq} by
\begin{align*}
  \rho \frac{\partial^2 u}{\partial t^2}
  +\alpha_u\frac{\partial u}{\partial t}-\DIV \left((1-z)^2 \sigma [u] \right) = f(t)\quad &\mbox{in}~\Omega\times [0,T],
\end{align*}
where $\rho\ge 0$ and $\alpha_u\ge 0$ are the material's density and friction coefficient,
respectively.
The coefficient $\alpha_{u}$ represents
the friction between the elastic body and a stationary background.
In a normal three-dimensional elastic body, friction with the background does not exist.
However, in a two-dimensional setting, contact friction can occur.
For example, when breaking an elastic sheet on a substrate (lubricated plate)
\cite{Hayakawa}, $\alpha_{u}$ has a positive value, which is controllable by the experiment.

When the friction is significant, and the inertia is negligible,
we can assume $\rho=0$ and $\alpha_u>0$. Then, the elliptic force balance equation
is replaced by the parabolic one. This model was proposed in
\cite{T-K09}, where $\alpha_u>0$ was considered a small coefficient
to regularize the elliptic equation. Since the elliptic linear elasticity equation
is degenerated if the damage variable $z$ has a value of $1$,
and it is numerically unstable even if the value of $z$ is very close to $1$,
the small parameter $\alpha_u >0$ is helpful to get a stable numerical solution.
However, $\alpha_u$ has an effect of a mathematical or numerical regularization
and a physical meaning as a coefficient of friction.

\section{Conclusion}\label{Sec:conclusion}
\setcounter{equation}{0}
We have shown that the time relaxation parameter $\alpha >0$
in the fracture phase field model has a concrete physical meaning as
$\alpha =\frac{d G_c^*}{dV}$, where $G_c^*$ is the velocity-dependent fracture energy
of the material, and $V$ is the crack tip velocity.
Such velocity dependence of the fracture energy is caused by a
process zone formation near the crack tip and is deeply related to the physical 
energy dissipation. In conclusion, the small parameter $\alpha >0$ in the F-PFM
is not only for the mathematical stabilization
of the variational fracture model
but also a physical quantity related to the energy dissipation process
during the crack propagation.

Furthermore, we derived the Griffith-type fracture model \eqref{aV}
for nonlinear $V$-dependence of the fracture energy: $G_c^*(V)=G_c+\alpha^* (V)$,
and proved the well-posedness of the model (Remark~\ref{wp}).
According to experimental measurements of some polymers, e.g. \cite{TFM2000},
the $V$-dependence is not always linear but exhibits
several nonlinearities.

In this study, we have established that the F-PFM
corresponds to the case of linear $V$-dependence: $G_c^*(V)=G_c+\alpha V$.
This analysis suggests a further generalization of the F-PFM with
a nonlinear $V$-dependent $G_c^*(V)=G_c+\alpha^*(V)$ 
\begin{align}\label{F-PFM2}
\alpha^* \left(\frac{\partial z}{\partial t}\right)=
\left( -\frac{\delta \cE_{tot}^*}{\delta z}\right)_+,
\end{align}
with \eqref{HGcV} as an analogy of \eqref{aV}.
The model \eqref{F-PFM2} is expected to be a potential mathematical model
for crack propagation in polymers and hydrogels, which often exhibit
nonlinear $V$-dependence of the fracture energy.

In conclusion, we revealed that the physical origin of the gradient flow structure
of the variational fracture models is the velocity dependence of
the fracture energy, which is originated from the localized energy dissipation 
by the formation of the process zone around the crack tip.

\appendix
\section{A lemma for the positive part}\label{Sec:lem-pp}
\setcounter{equation}{0}
Let us define the positive part of $c\in\R$ by
$(c)_+:=\max (c,0)$.
We repeatedly use the following simple lemma concerning the positive part
in our arguments. 
\begin{Lem}\label{lemabc}
For $a,~b\in\R$, it holds that
\begin{align}\label{ab+}
  a=(a-b)_+\quad \Longleftrightarrow \quad
  \left\{
  \begin{array}{l}
    a\ge 0\\ b\ge 0\\ ab=0,
    \end{array}\right.
\end{align}
Alternatively, for $a,~c\in\R$, it holds that
\begin{align}\label{ac+}
  a=(c)_+\quad \Longleftrightarrow \quad
  \left\{
  \begin{array}{l}
    a\ge 0\\ a\ge c\\ a(a-c)=0.
  \end{array}\right.
\end{align}
\end{Lem}
\proof
As the relation \eqref{ab+} is obtained from \eqref{ac+}
by replacing $b=a-c$, we prove \eqref{ac+}.
Suppose $a=(c)_+$. Then $a\ge 0$ and $a\ge c$ hold. 
If $a-c>0$, it implies $c<a=(c)_+$ and $a=(c)_+=0$.
Hence, the three conditions on the right-hand side of \eqref{ac+}
are derived.

Conversely,
if we suppose the three conditions on the right-hand side of \eqref{ac+},
one of the following two cases holds:
(i) $0<a=c$, (ii) $0=a\ge c$,
in both cases, we can quickly check that the condition $a=(c)_+$ holds.
\qed

\section{Crack evolution and crack path}\label{Sec:cecp}
\setcounter{equation}{0}

According to Section~2 of \cite{AKA}, we define admissible sets of cracks:
\begin{equation*}
\begin{aligned}
  \cC &:= \left\{ \Sigma \subset \Omega\ ;~\Omega\setminus\Sigma
  \text{ is open},~\overline{\Sigma} \cap \overline{\Gamma_D }
  = \varnothing,~\cH^{d-1}(\Sigma) < \infty\right\},\\
  \cC_0 &:= \left\{ \Sigma\in\cC\ ;\ \Omega\setminus\Sigma\text{ is connected}\right\},\\
\end{aligned}
\end{equation*}
where $\cH^{d-1}$ is the $(d-1)$-dimensional Hausdorff measure
and we set $|\Sigma|:=\cH^{d-1}(\Sigma)$.
\begin{Def}[Crack evolution]\label{def:crackEvolution}
  If $\{\Sigma(t)\}_{t_0\leq t\leq t_1}$ satisfies the following conditions,
  we call it a crack evolution in $\Omega$.
  1) $\Sigma(t) \in \cC_0$ $(t_0\leq t < t_1)$, $\Sigma(t_1)\in \cC$.
  2) $\Sigma(t)\subset\Sigma(\Tilde{t})$ $(t_0\leq t\leq \Tilde{t} \leq t_1)$.
  3) $\cH^{d-1}(\Sigma(t)) = \cH^{d-1}(\Sigma(\Tilde{t}))$
    implies $\Sigma(t) = \Sigma(\Tilde{t})$.
  Furthermore, if $\cH^{d-1}(\Sigma(t))$ is continuous within $t\in [t_0, t_1]$,
  then $\{\Sigma(t)\}_{t_0\leq t\leq t_1}$ is called a continuous crack evolution in $\Omega$.
\end{Def}

\begin{Def}[Crack path]\label{def:crackPath}
  If $\{\Sigma_p(\ell)\}_{\ell_0 \leq \ell \leq \ell_1}$ satisfies the following conditions,
  we call it a crack path in $\Omega$.
1) $\Sigma_p(\ell) \in \cC_0$ $(\ell_0 \leq \ell < \ell_1)$, $\Sigma_p(\ell_1) \in \cC$.
2) $\Sigma_p(\ell) \subset \Sigma_p(\Tilde{\ell})$
 $(\ell_0 \leq \ell \leq \Tilde{\ell} \leq \ell_1)$.
3) $\ell = \cH^{d-1}(\Sigma_p(\ell))$ $(\ell_0 \leq \ell \leq \ell_1)$.
\end{Def}

\begin{Prop}\label{prop:crackPathProp}
  If $\{{{\Sigma}(t)}\}_{t_0 \leq t \leq t_1}$ is a continuous crack evolution in $\Omega$,
  then there exists a unique crack path $\{\Sigma_p(\ell)\}_{\ell_0 \leq \ell \leq \ell_1}$
  in $\Omega$ such that $\Sigma(t)= \Sigma_p(L(t))$ for
  $t\in [t_0, t_1],$ where $L(t):=\cH^{d-1}(\Sigma(t))$.
\end{Prop}
There is a proof of this proposition in Section~2 of \cite{AKA}.

~\\\noindent
{\bf Acknowledgement}: 
This work was partially supported by JSPS KAKENHI Grant Nos. 20KK0058,
20H01812, and 21K03356.

\end{document}